\documentclass{etds}

\newtheorem{theo}{Theorem}[section]

\newtheorem{lemma}{Lemma}[section]
\newtheorem{cor}{Corollary}[section]

\begin{document}
\ETDS{1}{38}{?}{2002}

\runningheads{D. Sz\'asz and T. Varj\'u}{LLT for the Lorentz Process and Its Recurrence in the Plane}

\title{Local Limit Theorem for the Lorentz Process and Its Recurrence in the Plane}

\author{Domokos Sz\'asz, Tam\'as Varj\'u\footnote{Research supported by the Hungarian National Foundation for Scientific
    Research grants No. T32022 T26176 and Ts040719}}

\address{Budapest University of Technology\\   Mathematical Intitute
and Center of Applied Mathematics \\ Budapest, Egry J. u. 1 Hungary H-1111\\
\email{szasz@renyi.hu, kanya@math.bme.hu}}

\recd{3 May 2002}

\begin{abstract}
  For Young systems, i.\ e.\ for hyperbolic systems without/with singularities satisfying Young's axioms \cite{Young}
  (which imply exponential decay of correlation and the CLT) a local CLT is proven. In fact, a unified version of the local CLT
  is found, covering among others the absolutely contionuous and the arithmetic cases. For the planar Lorentz process with
  a finite horizon this result implies a.) the local CLT and b.) the recurrence.  For the latter case ($d=2$, finite
  horizon), combining the global CLT with abstract ergodic theoretic ideas, K. Schmidt, \cite{Sch} and J.-P. Conze,
  \cite{Conze}, could already establish recurrence.
\end{abstract}

\section{Introduction}
\subsection{Motivation}

The \emph{Lorentz process} is a physically utmost interesting mechanical model of Brownian motion (cf.\ \cite{Sz00}). It
is the deterministic motion of a point particle starting from a random phase point and undergoing specular reflections on
the boundaries of strictly convex scatterers. Throughout this paper we will only consider a $\mathbb Z^d$-periodic
configuration of scatterers. Once it had been established that the diffusion limit of the planar Lorentz process is,
indeed, the Wiener process (\cite{BS}, see also \cite{BSCH91}), the question of its 
{\it recurrence} was immediately raised by Ya. G. Sinai.  Here
recurrence means that the process almost surely returns to any fixed bounded domains of the configuration space.  In fact,
for Lorentz processes the exact analogue of P\'olya's theorem known for random walks is strongly expected. The first
positive result was obtained in \cite{KSz1}, where a slightly weaker form of recurrence was demonstrated: the process
almost surely returns infinitely often to a moderately (actually logarithmically) increasing sequence of domains. The
authors used a probabilistic method combined with the dynamical tools of Markov approximations.  The weaker form of the
recurrence was the consequence of the weaker form of their local limit theorem: they could only control the probabilities
that the Lorentz process $S_n$ in the moment of $n^\mathrm{th}$ collision falls into a sequence of moderately increasing
domain rather than into a domain of fixed size. These results, moreover, were restricted to the finite horizon case, i.\ 
e.\ to the case when there is no orbit without any collision.

A novel - and surprising - approach appeared in 1998-1999, when independently Schmidt \cite{Sch} and Conze \cite{Conze}
were, indeed, able to deduce the recurrence from the global central limit theorem (CLT) of \cite{BS} by adding (abstract)
ergodic theoretic ideas. Their approach seems to be essentially restricted to the finite horizon case and to $d=2$.  Our
main aim is to return to the probabilistic-dynamical approach and - still for the finite horizon case - we can first prove
a true local central limit theorem (LCLT) for the planar Lorentz process $S_n$.

\subsection{Statement of theorems}

As a matter of fact, beyond treating just the Lorentz process we are also able to obtain a LCLT in a much wider setup.
Namely our LCLT is valid whenever Young obtains a CLT. Her systems, called in our paper as 
{\it Young systems}, are introduced in
subsection \ref{subsec:You}. Roughly speaking, these are systems $(X, T, \nu)$
\begin{enumerate}
  \item whose every power is ergodic;
  \item which satisfy several technical assumptions well-known from hyperbolic theory;
  \item whose phase space $X$ contains a subset $\Lambda$ with a hyperbolic product structure;
  \item where the return time into $\Lambda$ has an exponentially decaying tail.
\end{enumerate}

For stating our main theorem we have to fix some notations first. For a fixed $f:X \to \mathbb{R}^d$ denote the average
$\nu (f) = a$, and \[ S_n(x) = \sum _{k=0}^ {n-1} f (T^{k}x) \] the Birkhoff sum. Consider the smallest translated closed
subgroup $V+r\subseteq\mathbb{R}^d$ which supports the values of $f$ ($V$ is the group and $r$ is the translation). By
ergodicity of all powers of $T$, the support of $S_n$ is $V+nr$.

\begin{theo}
  Suppose that
  \begin{enumerate}
    \item $(X, T, \nu)$ is a Young system (cf.\ subsection \ref{subsec:You});
    \item $f$ is minimal: i.\ e.\ it is not cohomologous to a function for which the support in the above sense is
      strictly smaller.
    \item $f$ is nondegenerate: i.\ e.\  $\text{span} \left<V\right>=\mathbb{R}^d$, and
    \item $f$ is bounded and H\"older-continuous.
  \end{enumerate}  
  Let $k_n\in V+nr$ be such that $\frac{k_n-na}{\sqrt{n}}\rightarrow k$.  Denote the distribution of $S_n-k_n$ by
  $\upsilon_n$. Then \[ n^\frac d2 \upsilon_n\rightarrow \varphi(k) l\] where $\varphi$ is a non-degenerate normal density
  function with zero expectation, and $l$ is the uniform measure on $V$: product of counting measures and Lebesgue
  measures. The convergence is meant in the weak topology.
\end{theo} 

\proc{Remark}
  For non-minimal functions we can obtain an analogous result. The limit measure on the right hand side in this case is
  not necessarily uniform.
\medbreak
\proc{Remark}
  Traditionally one formulates the LCLT for the absolutely continuous and for the arithmetic case separately. An advantage
  of our statement is that it is unified and beyond these two cases it also contains the mixed ones. Though for the
  absolutely continuous case it is slightly weaker than the LCLT for densities, nevertheless our variant, for instance, is
  still amply sufficient to treat recurrence properties.
\medbreak

Turning to the Lorentz process, let us denote by $(M, S^{\mathbb R}, \mu)$ a two-dimensional 
{\it dispersing billiard} dynamical
system with a finite horizon, the usual factor of the Lorentz process, where $\mu$ is the natural invariant probability
measure (the Liouville-one), and consider its Poincar\'e section $(\partial M, T, \mu_1)$ (for formal definitions of
billiards cf.\ section \ref{sec:rec}).

In case one takes $f$ as $\kappa: \partial M \to \mathbb R^2$, the discrete free flight function of the planar Lorentz
process, then this result combined with considerations of \cite{KSz1}, and an asymptotic independence statement proved
right after the main theorem immediately provide the recurrence of $S_n$ as well. It will be shown in section
\ref{sec:rec} that $\kappa$ satisfies the conditions of the main theorem.

\begin{cor}
  The planar Lorentz process with a finite horizon is almost surely recurrent.
\end{cor}

\subsection{Some history}

LCLT's for functions of a Markov chain were first obtained by Kolmogorov in 1955 using probabilistic ideas. Then, in 1957,
Nagaev, \cite{Nag} -- by using operator valued Fourier transforms and perturbation theory -- could find a general form of
LCLT's for functions of a Markov chain.  Independently, variants of this method got later rediscovered and/or applied A) by
Kr\'amli and Sz\'asz \cite{KSz2} to prove a LCLT for random walks with internal states, B) by
 Guivarch and Hardy \cite{GH}
in the setting of Anosov diffeomorphisms  C) by Roussean-Egele, \cite{RE83}, Morita \cite{M94} and Broise \cite{Bro} for
expanding maps of the interval and finally D) by Aaronson and Denker, \cite{Aa} in the setting of Gibbs-Markov maps.

Beyond establishing LCLT's for the planar Lorentz process for the first time, the technical interest and achievement of
this paper is the following: \cite{BS} and \cite{KSz1} used a Markov approximation scheme of the Lorentz process based
upon the Markov partition of the Sinai billiard. Several later works demonstrated that a Markov partition for a hyperbolic
system with singularities is a too rigid construction, and introduced Markov sieves \cite{BSCH91} and finally Markov
returns \cite{Young} instead. Our aim therefore is to work out {\it how Markov returns can be used to prove 
probabilistic statements}
 (e.\ g.\ to a large deviation result we return in a forthcoming paper).

\subsection{Probabilistic 
ideas: What is a local CLT and what is its relation to recurrence?}

For illustrating the probabilistic ideas, take a simple symmetric 
random walk (SSRW) on $\mathbb Z^d$.
 So let
 $W_n = X_1 + \dots X_n$, where $X_1, \dots, X_n, \dots$
are independent, identically distributed random variables with the common distribution
$P(X_i = \pm e_j) = \frac{1}{2d};\ 1 \le j \le d$ for all $i \in 
\mathbb Z_+$ (here the $e_j$s are the standard unit vectors of $\mathbb Z^d$).
 To investigate whether the SSRW is recurrent or not one turns to the Borel-Cantelli lemma: if 
$\sum_n P\left(W_n=0\right)=\infty$, then $P\left( \exists n_k\rightarrow \infty 
\ \
{\text {\rm such that}}\ \ \forall k\ W_{n_k}=0 \right) =1$. If the sum is convergent, then
 $P \left( \exists N 
\ \ {\text {\rm such that}}\ \ 
 \forall n>N \quad W_n\not=0\right)=1$. 
To apply the lemma we need to calculate the asymptotics of the probability $P(W_n=0)$. 

The CLT says that $P({W_n} \in {\sqrt n} A) \to \Phi (A)$ as $n \to \infty$, where $\Phi(A) = \int_A
\phi(s)ds$ and $\phi(s) = \exp{(-\frac{s^2}{2})}$ is the 
$d$-dimensional Gaussian density.
In other words it describes the asymptotics of a sequence of sets increasing like $\sqrt n$.
In contrast, for the SSRW  the local CLT  
says that,  $ n^{d/2} P(W_n = [s\sqrt n]) \to 2 \phi(s)$ as $n \to 
\infty$,  if the sum of the coordinates of the vector $[s \sqrt n] - n$ is even. In other words,
the LCLT describes the asymptotics of a sequence of sets of fixed size (in this case: a point),
consequently it is, indeed, local! As an application we get $P(W_n=0)\sim 2 n^{-d/2}$, 
consequently for the line and the plane SSRW is recurrent, for higher dimension it is transient.

This paper is organized as follows. Primarily, in section 2, we will formulate the abstract setting, define the notion of
Young-systems and recall our basic spectral tool: the Doeblin-Fortet (in the theory of dynamical systems also known as
Lasota-Yorke) inequality. Section 3 is devoted to important spectral properties of the Fourier transform:
quasicompactness, arithmeticity and a useful Nagaev-type theorem on a one-dimsnional approximation of the Fourier
transform in a neighbourhood of the origin.  In section 4 we establish our local limit theorem for Young systems and, in
addition, a certain asymptotic independence statement necessary to prove the recurrence.  In the fifth section we turn our
attention to billiards and to the Lorentz-process (with prerequisites in subsection 5.1) to get recurrence of planar
Lorentz-process as an application of the abstract theorems in subsection 5.3. In subsection 5.2 we analyze the
arithmeticity of the discrete free flight function.

\section{Prerequisites}

Since local central limit theorems are refined versions of (global) central limit theorems, it is not surprising that our approach
relies heavily on Young's work \cite{Young}, where -- among others -- an exponential decay of correlations and a central
limit theorem were proved for 2-$D$ dispersing billiards with a finite horizon. Here we present a concise summary of the
main points of Young's paper, which are necessary for our consideration.

\subsection{Young systems}\label{subsec:You}

Let $T$ be a $C^{1+\epsilon}$ diffeomorphism with singularities of a compact Riemannian manifold $X$ with boundary. More
precisely, there exists a finite or countably infinite number of pairwise disjoint open regions $\{X_i\}$ whose boundaries
are $C^1$ submanifolds of codimension 1, and finite volume such that $\cup X_i = X$, $T\big|_{\cup X_i}$ is $1-1$ and
$T\big|_{X_i}$ can be extended to a $C^{1+\epsilon}$-diffeomorphism of $\bar X_i$ onto its image. Then $\breve S = X
\setminus \cup X_i$ is the \emph{singularity set}. Later, for billiards, we will also use the notation $S = \breve S \cup
T^{-1}\breve S$. The Riemannian measure will be denoted by $\mu$, and if $W\subset X$ is a submanifold, then $\mu_W$ will
denote the induced measure. The invariant Borel probability measure will be denoted by $\nu$.  \proc{Definition} An
embedded disk $\gamma\subset X$ is called an \emph{unstable manifold} or an \emph{unstable disk} if $\forall x,y\in\gamma,
\enskip d(T^{-n}x,T^{-n}y)\rightarrow 0$ exponentially fast as $n\rightarrow\infty$; it is called a \emph{stable manifold}
or a \emph{stable disk} if $\forall x,y\in\gamma, \enskip d(T^nx,T^ny)\rightarrow 0$ exponentially fast as
$n\rightarrow\infty$. We say that $\Gamma^u=\{\gamma^u\}$ is a \emph{continuous family of $C^1$ unstable disks} if the
following hold:
\begin{itemize}
  \item $K^s$ is an arbitrary compact set; $D^u$ is the unit disk of some $\mathbb{R}^n$;
  \item $\Phi^u\colon K^s\times D^u\rightarrow X$ is a map with the property that
  \begin{itemize}
    \item $\Phi^u$ maps $K^s\times D^u$ homeomorphically onto its image,
    \item $x\rightarrow\Phi^u\mid(\{x\}\times D^u)$ is a continuous map from $K^s$ into the space of $C^1$ embeddings of
      $D^u$ into X,
    \item $\gamma^u$, the image of each $\{x\}\times D^u$, is an unstable disk.
  \end{itemize}
\end{itemize}
\emph{Continuous families of $C^1$ stable disks} are defined similarly.
\medbreak
\proc{Definition}
  We say that $\Lambda\subset X$ has a \emph{hyperbolic product structure} if there exist a continuous family of
  unstable disks $\Gamma^u=\{\gamma^u\}$ and a continuous family of stable disks $\Gamma^s=\{\gamma^s\}$ such that
  \begin{enumerate}
    \item[(i)] $\dim \gamma^u+\dim\gamma^s=\dim X$
    \item[(ii)] the $\gamma^u$-disks are transversal to the $\gamma^s$-disks with the angles between them bounded away
      from $0$;
    \item[(iii)] each $\gamma^u$-disk meets each $\gamma^s$-disk in exactly one point;
    \item[(iv)] $\Lambda=(\cup\gamma^u)\cap(\cup\gamma^s)$.
  \end{enumerate}
\medbreak
\proc{Definition}
  Suppose $\Lambda$ has a hyperbolic product structure. Let $\Gamma^u$ and $\Gamma^s$ be the defining families for
  $\Lambda$. A subset $\Lambda_0\subset\Lambda$ is called an \emph{$s$-subset} if $\Lambda_0$ also has a hyperbolic
  product structure and its defining families can be chosen to be $\Gamma^u$ and $\Gamma^s_0$ with
  $\Gamma^s_0\subset\Gamma^s$; \emph{$u$-subsets} are defined analogously. For $x\in\Lambda$, let $\gamma^u(x)$ denote
  the element of $\Gamma^u$ containing $x$.
\medbreak
In general a measurable bijection $M:(X_1,m_1)\rightarrow(X_2,m_2)$ between two finite measure spaces is called
\emph{nonsingular} if it maps sets of $m_1$-measure $0$ to sets of $m_2$-measure $0$. If $M$ is nonsingular, we define the
Jacobian of $M$ wrt $m_1$ and $m_2$, written $J_{m_1,m_2}(M)$ or simply $J(M)$, to be the Radon-Nikodym derivative
$\frac{d(M^{-1}_*m_2)}{dm_1}.$ To denote $J(T)$ wrt $\mu_{\gamma^u}$ we will use $\det DT^u.$

\proc{Definition}
  We call $(X,T,\nu)$ a \emph{Young system}, if the following Properties {\bf (P1)-(P8)} are true:
\medbreak

\begin{enumerate}
\item[\bf (P1)] There exists a $\Lambda\subset X$ with a hyperbolic product structure and with
  $\mu_\gamma\{ \gamma\cap\Lambda \}>0$ for every $\gamma\in\Gamma^u$.
\item[\bf (P2)] There is a countable number of disjoint $s$-subsets $\Lambda_1,\Lambda_2,\dots\subset\Lambda$ such that
  \begin{itemize}
  \item on each $\gamma^u$-disk $\mu_{\gamma^u}\{(\Lambda\setminus\cup\Lambda_i)\cap\gamma^u\}=0$;
  \item for each $i$, $\exists R_i\in\mathbb{Z}^+$ such that $T^{R_i}\Lambda_i$ is a $u$-subset of $\Lambda$;
  \item for each $n$ there are at most finitely many $i$'s with $R_i=n$;
  \item $\min R_i\geq$ some $ R_0$ depending only on $T$
  \end{itemize}
\item[\bf (P3)] For every pair $x,y\in\Lambda$, we have a notion of \emph{separation time} denoted by $s_0(x,y)$. If
  $s_0(x,y)=n$, then the orbits of $x$ and $y$ are thought of as being ``indistinguishable'' or ``together'' through their
  $n^\mathrm{th}$ iterates, while $T^{n+1}x$ and $T^{n+1}y$ are thought of as having been ``separated.''
  (This could mean that the points have moved a certain distance apart, or have landed on opposite sides of a
  discontinuity manifold, or that their derivatives have ceased to be comparable.) We assume:
 \begin{enumerate}
   \item[(i)] $s_0\geq 0$ and depends only on the $\gamma^s$-disks containing the two points;
   \item[(ii)] the number of ``distinguishable'' n-orbits starting from $\Lambda$ is finite for each $n$;
   \item[(iii)] for $x,y\in\Lambda_i, \enskip s_0(x,y)\geq R_i+s_0(T^{R_i}x,T^{R_i}y);$
 \end{enumerate}
\item[\bf (P4)] Contraction along $\gamma^s$ disks. There exist $C >0$ and $\alpha < 1$ such that for
  $y\in\gamma^s(x),\enskip d(T^nx, T^ny)\leq C\alpha^n\enskip\forall n\geq 0$.
\item[\bf (P5)] Backward contraction and distorsion along $\gamma^u$. For $y\in\gamma^u(x)$ and $0\leq k\leq
  n<s_0(x,y)$, we have
  \begin{enumerate}
  \item[(a)] $d(T^nx,T^ny)\leq C\alpha^{s_0(x,y)-n}$;
  \item[(b)] \[\log\prod_{i=k}^n\frac{\det DT^u(T^ix)}{\det DT^u(T^iy)}\leq C\alpha^{s_0(x,y)-n}.\]
  \end{enumerate}
\item[\bf (P6)] Convergence of $D(T^i|\gamma^u)$ and absolute continuity of $\Gamma^s$.
  \begin{enumerate}
  \item[(a)] for $y\in\gamma^s(x)$,\[\log\prod_{i=n}^\infty\frac{\det T^u(T^ix)}{\det T^u(T^iy)}\leq C\alpha^n\quad\forall
    n\geq 0.\]
  \item[(b)] for $\gamma,\gamma'\in\Gamma^u$, if $\Theta\colon\gamma\cap\Lambda\rightarrow\gamma'\cap\Lambda$ is defined
    by $\Theta(x)=\gamma^s(x)\cap\gamma'$, then $\Theta$ is absolutely continuous and
    \[\frac {d(\Theta_*^{-1}\mu_{\gamma'})} {d\mu_\gamma} (x) = \prod_{i=0}^\infty \frac {\det DT^u(T^ix)} {\det
    DT^u(T^i\Theta x)}.\] 
  \end{enumerate}
\item[\bf (P7)] $\exists C_0>0$ and $\theta_0<1$ such that for some
  $\gamma\in\Gamma^u$,\[\mu_\gamma\{x\in\gamma\cap\Lambda:R(x)>n\}\leq C_0\theta_0^n\quad \forall n\geq 0;\]
\item[\bf (P8)] $(T^n,\nu)$ is ergodic $\forall n\geq1$.
\end{enumerate}

\medbreak Now we will define the \emph{Markov extension}. Let $R:\ \Lambda \to \mathbb Z_+$ be the function which is
$R_i$ on $\Lambda_i$, and let
\[\Delta\stackrel{\mathrm{def}}{=}\{(x,l): x\in\Lambda;\enskip l=0,1,\dots,R(x)-1\}\] and define \[F(x,l)=\left\{
  \begin{array}{ll}
(x,l+1) & \text{if}\quad l+1<R(x)\\
(T^Rx,0) & \text{if}\quad l+1=R(x)
  \end{array}\right. \]
We will refer to $\Delta_l$ as the $l^\mathrm{th}$ level of the tower $\Delta$. Young also has a construction for $\tilde
\nu$, the SRB-measure of the extension, for which the pushforward is $\nu$, and $J(F)\equiv 1$ except on
$F^{-1}(\Delta_0)$.

On the tower a Markov partition $\mathcal{D}$ can be defined, with the following properties:
\begin{enumerate}
\item[(a)] $\mathcal{D}$ is a refinement of the partition $\Delta_l$. ($\mathcal{D}_l$ denotes $\mathcal{D}|\Delta_l$.)
\item[(b)] $\mathcal{D}_l$ has only a finite number of elements and
each one is the union of a collection of $\Lambda_i$'s;
\item[(c)] $\mathcal{D}_l$ is a refinement of $F\mathcal{D}_{l-1}$;
\item[(d)] if $x$ and $y$ belong to the same element of $\mathcal{D}_l$, then $s_0(F^{-l}x,F^{-l}y)\geq l$;
\item[(e)] if $R_i=R_j$ for some $i\not=j$, then $\Lambda_i$ and $\Lambda_j$ belong to different elements
  of $\mathcal{D}_{R_i-1}$.
\end{enumerate}
Let $\Delta^*_{l,j}=\Delta_{l,j}\cap F^{-1}(\Delta_0)$. We think of $\Delta_{l,j}\setminus\Delta^*_{l,j}$ as ``moving
upward'' under $F$, while $\Delta_{l,j}^*$ returns to the base.

\medbreak It is natural to \emph{redefine the separation time} to be $s(x,y)\stackrel{\mathrm{def}}{=}$ the
largest $n$ such that for all $i\leq n,\enskip F^ix$ and $F^iy$ lie in the same element of $\{\Delta_{l,j}\}$.  We claim
that \textbf{(P5)} is valid for $x,y\in\gamma^u\cap\Delta_{l,j}$ with $s$ in the place of $s_0$. To verify this, first
consider $x,y\in\Lambda$. We claim that $s(x,y)\leq s_0(x,y)$. If $x,y$ do not belong to the same $\Lambda_i$, then this
follows from rule (d) in the construction of $\mathcal{D}_l$; if $x,y\in\Lambda_i$, but $T^Rx,T^Ry$ are not contained in
the same $\Lambda_j$, then $s(x,y)=R_i+s(T^Rx,T^Ry)$, which is $\leq s_0(x,y)$ by property \textbf{(P3)},(iii) of $s_0$,
and so on. In general, for $x,y\in\Delta_{l,j}$, let $x_0=F^{-l}x,\enskip y_0=F^{-l}y$ be the unique inverse images of $x$
and $y$ in $\Delta_0$. Then by definition $s(x,y)=s(x_0,y_0)-l$, and what is said earlier on about $x_0$ and $y_0$ is
equally valid for $x$ and $y$.

\textit{From here on $s_0$ is replaced by $s$ and \textbf{(P5)} is modified accordingly.}

\medbreak Now we recall an important distorsion property of the so called sliding map.  Fix an arbitrary
$\hat{\gamma}\in\Gamma^u$. For $x\in\Lambda$, let $\hat{x}$ denote the point in $\gamma^s(x)\cap\hat{\gamma}$, and define
\[u_n(x)=\sum_{i=0}^{n-1}(\varphi(T^ix)-\varphi(T^i\hat{x}))\]where $\varphi = \log \left| \det DT^u \right|$. From
\textbf{(P6)}(a) it follows that $u_n$ converges uniformly to some function $u$. On each $\gamma\in\Gamma^u$, we let
$m_\gamma$ be the measure, whose density wrt $\mu_\gamma$ is $e^u\cdot 1_{\gamma\cap\Lambda}$. Clearly, $T^{R_i}|
(\Lambda_i\cap\gamma)$ is nonsingular wrt these reference measures. If $T^{R_i} (\Lambda_i\cap\gamma) \subset \gamma'$,
then for $x\in\Lambda_i\cap\gamma$ we write $J(T^R) (x) = J_{m_\gamma,m_ {\gamma'} } (T^{R_i}|(\Lambda\cap\gamma)) (x)$.

\begin{lemma}\label{lemma1}
  \begin{enumerate}
  \item[(1)] Let $\Theta_{\gamma,\gamma'}\colon\gamma\cap\Lambda\rightarrow\gamma'\cap\Lambda$ be the sliding map along
    $\Gamma^s$. Then $\Theta_*m_\gamma=m_{\gamma'}$.
  \item[(2)] $J(T^R)(x)=J(T^R)(y)\enskip\forall y\in\gamma^s(x).$
  \item[(3)] $\exists C_1>0$ such that $\forall i$ and $\forall x,y\in\Lambda_i\cap\gamma,$
    \[\left|\frac{J(T^R)(x)}{J(T^R)(y)}-1\right|\leq C_1\alpha^{\frac 12s(T^Rx,T^Ry)}.\]
  \end{enumerate}
\end{lemma}

Next Young uses a factorised dynamics with a factorisation along stable manifolds of $\Delta$.\@ The advantage is that
this dynamics will behave as an expanding map, a simpler object to study. Let $\bar{\Delta}:=\Delta/\sim$ where $x\sim y$
iff $y\in\gamma^s(x).$ Since $F$ takes $\gamma^s$-leaves to $\gamma^s$-leaves, the quotient dynamical system
$\bar{F}\colon\bar{\Delta}\rightarrow\bar{\Delta}$ is clearly well defined.

Let us define $\bar{m}$ in the following way: let $\bar{m}|\bar{\Delta}_l$ be the measure induced from the natural
identification of $\bar{\Delta}_l$ with a subset of $\bar{\Delta}_0$, so that $J(\bar{F})\equiv 1$ except on
$\bar{F}^{-1}(\bar{\Delta}_0),$ where $J(\bar{F})=J(\overline{T^R}\circ\bar{F}^{-(R-1)}).$

We now define $\bar{m}$ on $\bar{\Lambda}$ following the ideas that have been used for Axiom A.\@ Lemma \ref{lemma1}
(1) allows us to define $\bar{m}$ on $\bar{\Lambda}$ to be the measure whose representative on each
$\gamma\in\Gamma^u$ is $m_\gamma$. Statement (2) says that $J(T^R)$ is well defined wrt $\bar{m}$, and
(3) says that $\log J(T^R)$ has a dynamically defined H\"older type property, in the sense that
$\alpha^{s(T^Rx,T^Ry)}$ could be viewed as a notion of distance between $T^Rx$ and $T^Ry$ (see \textbf{(P5)}). By using
this lemma Young obtains a distorsion property of the factorised map with a weaker constant $\beta$. Let $\beta$ be such
that $\alpha^{\frac12}\leq\beta<1$, and let $C_1$ be as in Lemma \ref{lemma1} (3).
\begin{enumerate}
\item[(I)] Height of tower.
  \begin{enumerate}
  \item[(i)] $R\geq N$ for some $N$ satisfying $C_1e^{C_1}\beta^N\leq \frac1{100}$;
  \item[(ii)] $\bar{m}\{R\geq n\}\leq C_0'\theta_0^n\enskip \forall n\geq0$ for some $C_0'>0$ and $\theta_0<1$.
  \end{enumerate}
\item[(II)] Regularity of the Jacobian.
  \begin{enumerate}
  \item[(i)] $J\bar{F}\equiv1$ on $\bar{\Delta}-\bar{F}^{-1}(\bar{\Delta}_0)$,
  \item[(ii)]
    \[\left| \frac{J \bar{F} (\bar{x})} {J\bar{F} (\bar{y})} -1 \right| \leq C_1 \beta^ {s (\bar{F} \bar{x},\bar{F}
      \bar{y})} \quad \forall\bar{x},\bar{y}\in\bar{\Delta}^*_{l,j}.\]
  \end{enumerate}
\end{enumerate}
Young proves \cite{Young}, that there exists an invariant probability measure $\bar{\nu}$, absolutely continuous
wrt $\bar{m}$, such that $\rho=\frac {d\bar{\nu}} {d\bar{m}}$ is bounded away from zero and infinity, and is
Lipschitz-continuous wrt the distance $\beta^s$. 

\subsection{The Doeblin-Fortet inequality and spectral properties}
 
\proc{Definition}
  Let $(\mathcal{C},\mathcal{L})$ be a pair of Banach spaces, such that $\mathcal{L}\leq \mathcal{C}$ is a linear
  subspace, $\left\|\,.\,\right\|_{\mathcal{L}}\geq\left\|\,.\,\right\|_{\mathcal{C}}$. We call this pair \emph{adapted}
  if each $\mathcal{L}$-bounded set is precompact in $\mathcal{C}$.
\medbreak

\proc{Definition} Let $(\mathcal{C},\mathcal{L})$ be an adapted pair. We call an
$A\colon\mathcal{C}\rightarrow\mathcal{C}$ bounded linear operator a \emph{Doeblin-Fortet operator}, if
$\exists\tau<1,\exists K>0,\exists n\in\mathbb{N}\quad\forall \varphi\in\mathcal{L},$
  \[\left\|A^n\varphi\right\|_{\mathcal{L}}\leq \tau \left\| \varphi \right\|_{\mathcal{L}} + K \left\| \varphi
  \right\|_{\mathcal{C}}.\] This latter is called the \emph{Doeblin-Fortet inequality}.
\medbreak

\begin{theo}{\cite{IT-M}}\label{itm}
  If $A$ is a Doeblin-Fortet operator on the adapted pair $(\mathcal{C},\mathcal{L})$, then $\exists\vartheta<1,N\geq 1$,
  projections $E_1,\dots,E_N$ onto finite dimensional subspaces of $\mathcal{L}$, and $\lambda_1 ,\dots, \lambda_N \in \{
  z \in \mathbb{C} : \left| z \right| = 1 \}$ such that $\forall \varphi\in\mathcal{L},n\in\mathbb{N}$ \[ \left \| A^n
      \varphi - \sum_{k=1}^N \lambda_k^n E_k \varphi \right\|_{\mathcal{L}} \leq K \vartheta^n \left\| \varphi
    \right\|_{\mathcal{L}}.\]
\end{theo}

Now we will define the \emph{function spaces}.  Let $\epsilon>0$ be such that
\begin{enumerate}
\item[($\epsilon$i)] $e^{2\epsilon}\theta_0<1,$
\item[($\epsilon$ii)] $\bar{m}(\bar{\Delta}_0)^{-1}\sum_{l,j}\bar{m}(\bar{\Delta}^*_{l,j})e^{l\epsilon}\leq2.$
\end{enumerate}
Now we are ready to define the function spaces. The elements will be functions
$\bar{\varphi}\colon\bar{\Delta}\rightarrow\mathbb{C}$ and the $\mathcal{C}$ norm is\[
\left\|\bar{\varphi}\right\|_{\mathcal{C}}\stackrel{\mathrm{def}}{=}\sup_{l,j}\left\|\bar{\varphi}|_{\bar{\Delta}_{l,j}}\right\|_\infty
e^{-l\epsilon}\] where $\left\|\,.\,\right\|_\infty$ is the essential supremum wrt $\bar{m}$. By ($\epsilon$i) it is clear
that constant multiple of this norm dominates the $L_1$-norm wrt $\bar m$. Let us introduce
\[\left\|\bar{\varphi}\right\|_h\stackrel{\mathrm{def}}{=}\sup_{l,j}\left(
  \sup_{\bar{x},\bar{y}\in\bar{\Delta}_{l,j}}\frac{\left|\bar{\varphi}
(\bar{x})-\bar{\varphi}(\bar{y})\right|}{\beta^{s(\bar{x},\bar{y})}}\right) 
e^{-l\epsilon};\] where the inner $\sup$ is again essential supremum wrt 
$\bar{m}\times\bar{m}$  and $\mathcal{L}$-norm is
 \[\left\|\bar{\varphi}\right\|_{\mathcal{L}}\stackrel{\mathrm{def}}{=}\left\|\bar{\varphi}\right\|_{\mathcal{C}}+\left\|\bar{\varphi}\right\|_h.\]
 
 $\mathcal{C}$ resp.\ $\mathcal{L}$ consist of functions for which the $\mathcal{C}$-norm resp.\ $\mathcal{L}$-norm is
 finite. The adaptedness is an easy consequence of the Arzela-Ascoli theorem. The Perron-Frobenius operator acting on
 these spaces is defined as follows: \[P(\bar{\varphi}) (\bar{x}) = \sum_ {\bar{x}^{-1} : \bar{F}\bar{x}^{-1} = \bar{x}}
 \frac {\bar{\varphi} (\bar{x}^{-1})} {J\bar{F} (\bar{x}^{-1})}. \] This is the adjungate operator of $\bar \varphi
 \mapsto \bar \varphi \circ \bar F$ on $L_2(\bar m)$. By ($\epsilon$i) both $\mathcal{C}$ and $\mathcal{L}$ is contained
 in $L_2(\bar m)$. The fact, that $P$ is a bounded operator on $\mathcal{C}$ follows from ($\epsilon$ii). The similar
 statement for $\mathcal{L}$ is proved in \cite{Young}, where Young deduces that
\begin{enumerate}
\item[(i)] $P$ is a contraction in $\mathcal{L}$.
\item[(ii)] it satisfies the D-F inequality,
\item[(iii)] by Theorem \ref{itm} it has a spectral gap,
\item[(iv)] and by {\bf(P8)} its only eigenvalue on the unit circle is $1$ and it is simple. (The eigenfunction is the
  invariant density $\rho$.)
\end{enumerate}
Later we will need the adjungate of $\bar \varphi \mapsto \bar \varphi \circ \bar F$ on $L_2(\bar \nu)$, this is $P^\rho
(\bar \varphi) \stackrel {\mathrm{def}} = \frac 1 \rho P(\rho \bar \varphi)$. Note that the spectrum of $P$ and $P^\rho$
is the same, just the eigenfunctions are divided by $\rho$.

\section{Spectral properties of the Fourier-transform}

In this section we are working with Young systems throughout. Let $f\colon X\rightarrow\mathbb{R}^d$ be a bounded,
piecewise $\eta$-H\"older function i.\ e.\ $f(x)-f(y) \leq C_f d(x,y)^\eta$ whenever $x,y\in X_i$. We are going to
associate a function $\bar f \colon \bar \Delta \rightarrow \mathbb{R}^d$ of the symbolic space. First we \emph{pull back}
$f$ along the projection map $\pi \colon \Delta \rightarrow \cup T^n\Lambda$ to a function $\tilde{f} \colon \Delta
\rightarrow \mathbb{R}^d$. This is clearly bounded and by \textbf{(P5)} $\tilde f (x) - \tilde f (y) \leq C_f \left( C
  \alpha^{s(x,y)} \right)^\eta$ meaning $\tilde f$ is $\eta$-H\"older wrt the metric $\alpha^s$. Next we use a standard
method described for example in \cite{PP}. We choose an unstable manifold in each Markov-rectangle $\Delta_{l,j}$, and
consider the projection $\Xi$ which sends each point along its stable manifold to our preferred unstable manifold.
Consider the function $h \stackrel {\mathrm{def}}{=} \sum_{n=0}^\infty \left( \tilde f \circ F^n - \tilde f \circ F^n
  \circ \Xi \right)$! The defining series converges since $\tilde f F^nx - \tilde f F^n\Xi x \leq C_f d(T^n \pi x, T^n \pi
\Xi x)^\eta$ and by {\bf(P4)} $\leq C_f (C \alpha^n)^\eta$.
\begin{gather*}
  h - h \circ F = \sum_{n=0}^\infty \left( \tilde f \circ F^n - \tilde f \circ F^n \circ \Xi \right) -
\sum_{n=0}^\infty \left( \tilde f \circ F^n+1 - \tilde f \circ F^n \circ \Xi \circ F \right) \\ = \tilde f - \left[ \tilde
  f\circ\Xi +\sum_{n=0}^\infty \tilde f\circ F^{n+1}\circ\Xi - \tilde f\circ F^n\circ\Xi\circ F \right] .
\end{gather*}
This can be rewritten as $h - h \circ F = \tilde f - \bar f$, where $\bar f$ is defined by the expression in square
brackets. Evidently $\bar f$ is constant when restrticted to any stable manifold, so it can be regarded as a function
defined on $\bar \Delta$.
  \begin{lemma}\label{contmod}
    If $f\colon X\rightarrow \mathbb{R}^d$ is piecewise $\eta$-H\"older, and
    $\beta$ satisfies $1 >\beta\geq\alpha^{\eta/2}$, then the associated function
    $\bar{f}\colon \bar{\Delta}\rightarrow\mathbb{R}^d$ is bounded and Lipschitz-continuous wrt the metric $\beta^s$:
    \[\left|\bar{f}(\bar{x})-\bar{f}(\bar{y})\right|\leq C\beta^{s(\bar x,\bar y)}.\]
  \end{lemma}
  \proc{Proof} Let $\bar x, \bar y \in \bar \Delta$ such that $s(x,y) \geq 2n$ then {\bf(P5)} ensures\[ \left|\tilde f F^k
    x - \tilde f F^k y \right|, \left|\tilde f F^k \Xi x - \tilde f F^k \Xi y \right| \leq C_f (C\alpha^{2n-k})^\eta,
  \quad 0 \leq k \leq n.\] For all $k>0$ {\bf(P4)} gives \[ \left|\tilde f F^nx - \tilde f F^n\Xi x \right|, \left| \tilde
    f F^ny - \tilde f F^n\Xi y \right| \leq C_f (C \alpha^n)^\eta. \] Hence $|h(x)-h(y)| \leq 2 C_f \sum_{k=0}^n
    (C\alpha^{2n-k})^\eta + 2C_f \sum_{k=n+1}^\infty (C \alpha^n)^\eta \leq \mathrm{const} C_f \alpha^{n\eta}$ given $1
    >\beta\geq\alpha^{\eta/2}$ the latter estimate $\leq \bar C_f \beta^s$. \ep \medbreak

\subsection{Quasicompactness}
\label{sec:qcft}

The purpose of this subsection is to prove the Doeblin-Fortet inequality for the Fourier transform of the Perron-Frobenius
operator: \[P_t(\bar{\varphi}):=P(e^{it\bar{f}} \bar{\varphi})\quad (\bar{\varphi}\in\mathcal{C})\] where
$f\colon X\rightarrow\mathbb{R}^d$ measurable, and $t\in\mathbb{R}^d$. Simpifying the notations for a fixed $t$ denote
$\omega=e^{i\left<t,\bar{f}\right>}$, so $P_t(\bar{\varphi})=P(\omega \bar{\varphi})$. For to prove the inequality we need
the assumption of H\"older continuity for the measurable $f$.
  \begin{lemma}
    If $f$ and $\beta$ satisfies the conditions of the previous lemma, then the operator $P_t$ satisfies the
    Doeblin-Fortet inequality $\forall t\in\mathbb{R}^d$.
  \end{lemma}
\proc{Proof}
  \[ \|P_t^n\bar{\varphi}\|_\mathcal{L} = \|P_t^n\bar{\varphi}\|_\mathcal{C} + \|P_t^n\bar{\varphi}\|_h \leq
  \|P_t^n\|_\mathcal{C} \|\bar{\varphi}\|_\mathcal{C} + \|P_t^n\bar{\varphi}\|_h.\] By ($\epsilon$ii)
  $\|P_t^n\|_\mathcal{C} \leq 2$, so we only have to bound the continuity modulus.
  \begin{gather*}
    P^n_t(\bar{\varphi})=P^n(\omega_n\bar{\varphi}) \text{ where }
    \quad\omega_n(\bar{x}):=\prod_{k=0}^{n-1}\omega(\bar{F}^k\bar{x}). 
    \intertext{It follows that} P_t^n(\bar{\varphi})(\bar{x})=\sum_{\bar{x}^{-n}:T^n\bar{x}^{-n}=\bar{x}}
    \frac{\omega_n(\bar{x}^{-n})\bar{\varphi} (\bar{x}^{-n})} {J\bar{F}^n(\bar{x}^{-n})}
  \end{gather*}
 
  If $\bar{x}$ and $\bar{y}$ lie in the same element $\bar{\Delta}_{l,j}$, then the inverse images can be coupled:
  $\bar{x}^{-n}_i$ and $\bar{y}^{-n}_j$ form a pair if $\forall k,0\leq k \leq n\quad \bar{F}^k(\bar{x}^{-n}_i)$ and
  $\bar{F}^k(\bar{y}^{-n}_j)$ belong to the same element of the Markov partition $\{\bar{\Delta}_{l,j}\}$. That this is
  really a coupling is ensured by (e) in the definition of $\mathcal{D}$. For notational simplicity suppose that the
  inverse images are numbered according to the coupling. We have then the following expression for the continuity modulus:
    \begin{gather*}
      \left| P_t^n (\bar{\varphi}) (\bar{x}) -P_t^n (\bar{\varphi}) (\bar{y}) \right| = \left| \sum_{\bar{x}^{-n}_i :
          \bar{F}^n\bar{x}^{-n}_i =\bar{x}} \frac{\omega_n (\bar{x}^{-n}_i) \bar{\varphi} (\bar{x}^{-n}_i)} {J\bar{F}^n
          (\bar{x}^{-n}_i)} -\frac {\omega_n(\bar{y}^{-n}_i) \bar{\varphi} (\bar{y}^{-n}_i)} {J\bar{F}^n (\bar{y}^{-n}_i)}
      \right|. \intertext{The right hand side can be written as $\left| I+II \right|$ where} I=\sum _{\bar{x}^{-n}_i :
        \bar{F}^n\bar{x}^{-n}_i =\bar{x}} \omega_n (\bar{x}^{-n}_i) \left( \frac{\bar{\varphi} (\bar{x}^{-n}_i)}
        {J\bar{F}^n (\bar{x}^{-n}_i)} -\frac { \bar{\varphi} (\bar{y}^{-n}_i)} {J\bar{F}^n (\bar{y}^{-n}_i)} \right), 
      \intertext{and} II=\sum _{\bar{x}^{-n}_i : \bar{F}^n\bar{x}^{-n}_i =\bar{x}} \frac {\bar{\varphi}
        (\bar{y}^{-n}_i)}{J\bar{F}^n (\bar{y}^{-n}_i)} \left( \omega_n (\bar{x}^{-n}_i) -\omega_n (\bar{y}^{-n}_i)
      \right). \intertext{The first quantity can be estimated as follows:} \left| I \right| \leq \sum _{\bar{x}^{-n}_i :
        \bar{F}^n\bar{x}^{-n}_i =\bar{x}} \left| \frac{\bar{\varphi} (\bar{x}^{-n}_i)} {J\bar{F}^n (\bar{x}^{-n}_i)}
        -\frac { \bar{\varphi} (\bar{y}^{-n}_i)} {J\bar{F}^n (\bar{y}^{-n}_i)} \right|
    \end{gather*}
    Young \cite{Young} gets her D-F inequality by estimating the same quantity in the case where $n=N$. For the estimate
    of the second term we have to say something about the continuity modulus of $\omega$: 
    \begin{align*}
      \left| \omega(a) - \omega(b) \right| = \left| e^{i\left<t,\bar{f}(a)\right>} - e^{i\left<t,\bar{f}(b)\right>}
      \right| &\leq \left| t \right| \left| \bar{f}(a) - \bar{f}(b) \right|.\intertext{By lemma \ref{contmod} this
        latter is} &\leq \left| t \right| C \beta^ {s(a,b)}.\intertext{Then the continuity modulus of $\omega_N$:}
      \left| \omega_N (\bar{x}^{-N}_i) - \omega_N (\bar{y}^{-N}_i) \right| &= \sum_{k=0}^{N-1} \left| \omega (\bar{F}^k
        (\bar{x}^{-N}_i)) - \omega (\bar{F}^k (\bar{y}^{-N}_i)) \right| \\
      &\leq \sum_{k=0}^{N-1} \left| t \right| C \beta^{s(\bar{F}^k (\bar{x}^{-N}_i),\bar{F}^k (\bar{y}^{-N}_i))} \\
      &= \sum_{k=0}^{N-1} \left| t \right| C \beta^{s(\bar{x},\bar{y})+N-k} \\
      &\leq \frac {\beta \left|t\right| C\beta^{s(\bar{x},\bar{y})}}{1-\beta}
    \end{align*}
    $II$ can be estimated by taking absolute value term by term. Then the continuity modulus is multiplied by $P^N
    |\varphi| y\leq e^{\epsilon l}\left\|P^N|\varphi|\right\|_\mathcal{C}\leq
    e^{l\epsilon}2^N\left\|\varphi\right\|_\mathcal{C}$. From these it is easy to see, that in the D-F inequality this
    estimate of $II$ contributes to the coefficient of $\|\varphi\|_\mathcal{C}$ by $\frac {2^N \beta C |t|} {1-\beta}$,
    so it doesn't bother Young's estimate of $I$.  \ep \medbreak

\subsection{Minimality}
\label{sec:min}

Next we have to investigate the $t$ values, for which $P_t$ has an eigenvalue on the unit circle. Othervise $P_t$ is
strictly contractive by quasicompactness. As we will see, this 
is the question of minimality. First we give the basic definitions for an
arbitrary dynamical system $(X,T,\nu)$. Then we will investigate Young's symbolic system $(\bar\Delta,\bar F,\bar\nu)$ to
get a characterisation of the abovementioned $t$-values. Finally, we will prove that the definitions for a Young system
$(X,T,\nu)$, and for the associated symbolic system $(\bar \Delta,\bar F, \bar \nu)$ provide  the same answer. Thus 
we can characterise
the ``bad'' $t$-values, by concentrating on the minimality of our function on the original system.
\proc{Definition}
  We say that $f$ is \emph{cohomologous} to $g$ (notation: $f\sim g$) if $\exists h $ measurable such that $f-g=h-h\circ
  T$. Under the \emph{minimal support} of a function 
$f$ (notation: $S(f)$) we mean the   minimal translated closed subgroup of
  $\mathbb{R}^d$, which supports its values. We call a  translated closed subgroup the \emph{minimal lattice} of $f$ if
  it is the intersection of minimal supports in the cohomology class of $f$ ($M(f)=\cap_{g:g\sim f} S(g)$). We call $f$
  \emph{minimal} if $S(f)=M(f)$. We call $f$ \emph{degenerate}
 if $M(f)$ is contained in a smaller dimensional affine subspace
  of $\mathbb{R}^d$.
\medbreak
\begin{lemma}
  Fix the function $f$. Then $P^\rho_t\bar g=\lambda \bar g$ with $|\lambda|=1 \iff e^{it\bar f} \bar g=\lambda \bar
  g\circ \bar F$. Moreover $\bar g$ can be supposed to take values on the unit circle.
\end{lemma}
\proc{Proof}
  \begin{itemize}
  \item [$\Longrightarrow$] If $P^\rho_t\bar g=\lambda \bar g$ then by ($\epsilon$i) $\bar g\in\mathcal{L}\Longrightarrow
    \bar g\in L_2(\bar{m})$, and also $\bar{g}\in L_2(\bar{\nu})$ we can take: \[ \left< e^{it\bar{f}} \bar g, \bar g\circ
    \bar{F} \right>_{\bar\nu} = \left< P\left(e^{it\bar{f}} \bar g\right), \bar g \right>_{\bar\nu} =\left<\lambda \bar g,
    \bar g \right>_{\bar\nu} = \lambda \left\| \bar g \right\|_{L_2(\bar\nu)}^2. \] From Cauchy-Schwartz inequality it
    follows that $e^{it\bar{f}} \bar g=\lambda \bar g\circ \bar{F}$.  By ergodicity we can suppose $|\bar g|\equiv 1$.
  \item [$\Longleftarrow$] If $e^{it\bar{f}} \bar g=\lambda \bar g\circ \bar{F}$ then $P^\rho_t(\bar g) = \frac 1 \rho
    P(\rho e^{it\bar f}\bar g) = \frac \lambda \rho P(\bar g\circ \bar F \rho)= \lambda \bar g \frac {P(\rho)} \rho =
    \lambda \bar g$. Since $|\bar g|=1\Longrightarrow \bar g\in \mathcal{C}$, then it follows that $\bar g\in\mathcal{L}$
    \cite{IT-M}.
  \end{itemize}
  \ep \medbreak This lemma shows that the $t$ values for which the abovementioned property holds form a closed subgroup of
  $\mathbb{R}^d$, moreover the eigenvalues and -functions preserve the group structure. If $P^\rho_{t_1}\bar g_1 =
  \lambda_1 \bar g_1 \circ \bar F$ and $P^\rho_{t_2} \bar g_2 = \lambda_2 \bar g_2 \circ \bar F$, then $P^\rho_ {t_1 +
    t_2} \bar g_1 \bar g_2 = \lambda_1 \lambda_2 (\bar g_1 \bar g_2) \circ \bar F$. Also, for $t \in G$, $t\mapsto \bar
  g_t$ and $t\mapsto \lambda_t$ are uniquely determined by ergodicity. (Here $G$ denotes the subgroup of $\mathbb{R}^d$
  formed by these $t$ values.) This uniqueness can be easily derived from the multiplicative structure, and the already
  known spectral picture for $P=P_0$. Since $\lambda_t$ is a multiplicative functional of $t$, so the logarithm is a
  linear one, and therefore $-i\log \lambda_t=tr$ for some $r$ real vector. (Taking the adequate branch of the logarithm.)
\begin{theo}\label{thm:min}
  $M(\bar f)=\widehat {\mathbb{R}^d/G}+r$. There exist minimal functions in each cohomology class. The minimal function is
  unique iff it is constant.
\end{theo}
\proc{Proof}
  \begin{itemize}
  \item [$\subset$] We are going to prove that $\forall t\in G, \forall x\in M(\bar f)\enskip$ one has 
$ e^{itx}=e^{itr}$.  Since
    $t\in G$ we have $e^{it\bar f}\bar g=\lambda \bar g\circ \bar F$. Taking the logarithm
    \begin{equation}
      \label{eq:gh}
      t\bar f \equiv -i\log \lambda +i\log \bar g -i\log \bar g\circ \bar F \pmod {2\pi}.
    \end{equation}
    Remember that the first term on the right hand side is $tr$. By denoting $h=i\log \bar g$ we get that $t\bar
    f-(Z+tr)=h-h\circ \bar F$ for some $Z$, which takes values in $2\pi\mathbb{Z}$. To lift it to vector valued equation
    let us denote $\vec h=\frac{th}{|t|^2}, \enskip \vec Z=\frac{tZ}{|t|^2}+\bar f^{t^\perp}-r^{t^\perp}$, we get that
    $\bar f\sim \vec Z+r$, and the right hand side takes values in $H=t^\perp\oplus\frac{2\pi t}{|t|^2}\mathbb{Z}+r$. By
    definition $H\supset M(f)$, and since $\forall x\in H \enskip e^{itx}=e^{itr}$ this is true for $\forall x\in M(\bar
    f)$.
  \item [$\supset$] We are going to prove that if for $t\in\mathbb{R}^d \enskip$ and $ \forall x\in M(\bar f) \enskip$ we
    have $ e^{itx}=e^{itr}$, then $t\in G$. The condition means that $\exists Z, \enskip Z\sim \bar f, \enskip S(Z)\subset
    t^\perp \oplus \frac{2\pi t}{|t|^2}\mathbb{Z} + r$. Combining the condition with the cohomological equation we get
    $e^{itZ}=e^{itr}=e^{it(\bar f-h+h\circ \bar F)}$. After rearranging one obtains $e^{it\bar
      f}e^{-ith}=e^{itr}e^{-ith\circ \bar F}$, and by the previous lemma $t\in G$.
  \item [$\exists$] Let us revisit the congruence (\ref{eq:gh}). Observe that $i\log \bar g$ is also a linear functional
    of $t$, so $i\log \bar g=ts$ for some $s:\bar{\Delta}\rightarrow\mathbb{R}^d$. The function $Z$ derived from this
    congruence is also linear in $t$, so $Z=tz$. Denote by $H$ the orthocomplement of the linear subspace generated by
    $G$. Recalling the definition of $r,s$ and $z$ we can see, that $r^H,s^H$ and $z^H$ can be arbitrary, so let the
    latter one agree with $\bar{f}^H$, and the others be $0$.  We get $\bar f-(z+r)=s-s\circ \bar F$. Consider now
    $S(z+r)$. In the definition of $Z$ we said that it takes values in $2\pi \mathbb{Z}$, but $Z=tz$ gives $\forall t\in G
    \enskip e^{it(z+r)}=e^{itr}$, so from the already proven part of the theorem it follows that $S(z+r)=M(\bar f)$.
    Uniqueness is obvious: if $M(\bar f)$ is not a single point, then taking any $h:X\rightarrow M(\bar f)$ nonconstant
    $\bar f-h+h\circ \bar F$ is also a minimal function, and by ergodicity is not equal to $\bar f$.

  \end{itemize}
\ep \medbreak

Let us remark, that $M(S_n)=M(\bar f)+(n-1)r$. One of the inclusions ($\subset$) is trivial, the other ($\supset$) follows
from ergodicity of iterates. Now we turn our attention to the point, that neither the Markov extension, nor the
factorisation changes minimality properties.

\begin{theo}\label{thm:minpr}
  If $g$ is minimal in the class of $f$, then $\bar g$ is minimal in the class of $\bar f$. $M(f)=M(\bar f)$.
\end{theo}
\proc{Proof} First we prove, that $f\sim g \Longrightarrow \tilde f\sim \tilde g$. Indeed $f-g=h-h\circ T\Longrightarrow
\tilde f - \tilde g = \tilde h-\tilde h\circ F$. Consider now the construction of $\bar f$. It is clear that this
construction preserves addition, so it is enough to prove, that if $f$ is null-cohomologous, then $\bar f$ also. This
means $f=h-h\circ F$. Putting this in the definition of $\bar f$, we see, that $\bar f=h\circ\Xi - h\circ\Xi\circ F$. The
function $h\circ\Xi$ is constant along stable lines, so $\bar f$ is null-cohomologous in the factorised system. So far we
have reached $f\sim g\Longrightarrow \bar f\sim \bar g$. If $g$ is minimal, then $S(g)=S(\bar g)$. From the formula
$\supset$ is trivial, and it cannot be strictly smaller, since $\tilde g\sim \bar g$ would contradict the minimality of
$g$. The only thing remained is to show that $S(\bar g)=M(\bar g)$. If not there would be an other function with smaller
support in the same class. Consider this on $\Delta$, as a function constant along stable manifolds! It would be in the
class of $\tilde g$, which would again contradict the minimality.  \ep \medbreak

\subsection{A Nagaev type theorem}
\label{sec:nag}

Expand now $P_t$ in a Taylor series around $t=0$! $P_t (\bar{\varphi}) = P(e^ {i\left<t , \bar{f} \right>} \bar{\varphi})
= P (\bar{\varphi}) + it P \left( \bar{f} \bar{\varphi} \right) - \frac{t^2}2 P \left( \bar{f}^2 \bar{\varphi} \right) +
o(t^2) \left\| \bar{f}^2 \bar{\varphi} \right\|_\mathcal{L}$.  From lemma \ref{contmod} it follows that the norm exists,
so the second order Taylor-expansion at zero makes sense. Let us denote the operator $\bar{\varphi} \mapsto P (\bar{f}
\bar{\varphi})$ by $M$ (mean) and $\bar{\varphi} \mapsto P (\bar{f}^2 \bar{\varphi})$ by $\Sigma$ (covariance).

Denote by $\lambda_t$ the leading -also simple- eigenvalue of $P_t$, (we know that $\lambda_0=1$) and by $\tau_t$ the
projection operator corresponding to $\lambda_t$. The invariant density $\rho$ is known to be bounded away from zero and
infinity, and is H\"older. We know that $\tau_0=\rho\bar m$, since $\rho$ is the invariant density.  Consider the second
order Taylor polynomial of these two objects:
\begin{align*}
  \lambda_t &= 1+ iat- b\frac{t^2}2 + o(t^2) \\
  \tau_t &= \rho\bar{m} + \eta t + \chi t^2 +o(t^2)
\end{align*}

By definition $\tau_t P_t = \lambda_t \tau_t$. Expressing the terms by the above equations and considering the
coefficients of $t$ and $t^2$ we get the following:
\begin{align*}
  i \rho\bar{m} M + \eta P &= \eta + ia \rho\bar{m} \\
  -\frac12\rho\bar{m} \Sigma + i \eta M + \chi P &= \chi + ia\eta - \frac{b\rho\bar{m}}2
\end{align*}
evaluating these on $\rho$ we get from the first that $a= \bar{m} M(\rho)$. We are allowed to suppose that $M(\rho)$ is a
constant. This is because if we change $\bar{f}$ to a cohomologous $\bar f'$ the maximal eigenvalue does not change. Just
like in the case of $P_t^\rho$ we will study a conjugated operator with the same spectrum. Let us solve the equation:
$P(\bar f\rho)-\int f d\nu= Pu-u$. This is solvable since the left hand side $\in\ker \bar m$. Let us consider $\bar
f'=\bar f-\frac u\rho+\frac{u\circ\bar F}{\rho\circ \bar F}$! This is clearly cohomologous to $\bar f$. Let us consider
$M'(\rho)=P(\bar f'\rho)=P(\bar f\rho)-Pu+P(\frac{u\circ\bar F}{\rho\circ \bar F}\rho)$. This latter term is $\frac u\rho
P\rho=u$. So by the definition of $u$ $M'(\rho)=\int f d\nu$ constant. Evaluating the second equation on $\rho$ we get
$b=\bar{m}\Sigma'(\rho)=\int \bar {f'}^2 d\bar \nu$, remember, that $a$ was the average of the function, now $b$ is some
second moment, and we can define covariance by $\sigma^2=b-a^2$. It is also remarkable, that $\sigma$ is the second
central moment of a function cohomologous to $\bar f$. If $f$ is nondegenerate, each such quadratic form (and consequently
$\sigma$) is nondegenerate also. We have proved the following theorem:
\begin{theo}\label{thm:nag}
  There are constants $\epsilon >0,$ $ K>0$ and $\theta<1$ and a function
  $\rho:(-\epsilon,\epsilon)^d\rightarrow\mathcal{L}$ such that \[\left\| P_t^nh -\lambda_t^n\rho_t\int_{\bar\Delta}
    hd\bar m \right\|_{\mathcal{L}}\leq K\theta^n\left\|h\right\|_{\mathcal{L}} \quad \forall \left| t \right| <\epsilon,
  \enskip n\geq 1,\enskip h\in \mathcal{L}, \] and $\rho_0=\rho$, $\lambda_t=1+ait-(\sigma^2+a^2) \frac {t^2}2 + o(t^2)$ .
\end{theo}

\section{Proof of the Main Theorem}
\label{sec:pr}

In this section we are still going to consider Young systems, in general. Without loss of generality (by adding a scalar)
we can suppose, that $r=0$, which means, that $M(f)$ is a closed subgroup of $\mathbb{R}^d$. It also means, that
$P_{t+u}=P_t$ if $u\in G$, so the $t$ values are actually taken from $\widehat{M(f)}=\mathbb{R}^d/G$. Later we will
concentrate on compact parts of this group.

\begin{lemma}[\cite{Aa}]\label{lem:cdf}
  Suppose that $\mathcal{K}$ is a compact set of $\mathcal{L}$ operators such that each element of $\mathcal{K}$ is a
  Doeblin-Fortet operator, and none of them has an $\mathcal{L}$-eigenvalue on the unit circle.Then $\exists K >0$ and
  $\theta<1$ such that \[ \left\| Q^n \right\|_\mathcal{L} \leq K\theta^n \quad \forall n\geq 1, \quad Q\in \mathcal{K}.
  \]
\end{lemma}

For to apply this lemma we have to cut out a neighborhood of zero. In it, however, theorem \ref{thm:nag} holds. Now we are
able to prove our main theorem.

\begin{theo}\label{thm:main}
  Suppose that
  \begin{enumerate}
    \item $(X, T, \nu)$ is a Young system (cf.\ subsection \ref{subsec:You});
    \item $f$ is minimal (cf.\ subsection \ref{sec:min});
    \item $f$ is nondegenerate  (cf.\ subsection \ref{sec:min});
    \item $f$ is bounded and H\"older-continuous.
  \end{enumerate}  
  Let $k_n\in M(f)$ be such that $\frac{k_n - na} {\sqrt{n}}\rightarrow k\in\mathbb{R}^d$. Denote the distribution of
  $S_n-k_n$ by $\upsilon_n$, then \[ \lim_{n\rightarrow\infty} n^{\frac d2} \upsilon_n
  = \frac {e^{\frac {-k^2}{2\sigma^2}}}{\det\sigma\sqrt{(2\pi)}^d} l. \] $l$ is  the uniform
  distribution on $M(f)$, more exactly it  is product of suitable counting measures and Lebesgue measures.
\end{theo}

\proc{Proof}
  Suppose that we choose a random point of $X$ according to the invariant distribution $\nu$. Let the joint distribution
  of $(x,T^nx,S_n(x)-k_n)$ be denoted by $\Upsilon_n$! We are going to prove, that \[ \lim_{n\rightarrow\infty} n^{\frac
    d2} \Upsilon_n \rightarrow \frac {e^{\frac {-k^2}{2\sigma^2}}}{\det\sigma\sqrt{(2\pi)}^d} \nu^2\times l. \] The
  definitions of $\tilde\Upsilon_n$, and $\bar\Upsilon_n$ are straightforward. In the following paragraph we are going to
  see that it is enough to prove that the limit of $n^\frac d2\bar\Upsilon_n$ is $\bar\nu^2\times l$ multiplied by the
  gaussian density of covariance $\sigma$ at $k$.
  
  It is clear that a similar limit for $\tilde\Upsilon_n$ is sufficient. Remember the definition of $\bar f$! Nota bene if
  $f$ is minimal then $h$ takes values in $M(f)$. So in the language of $\Upsilon$ the
  factorisation means the application of the mapping $(x,y,\xi)\mapsto(x,y,\xi+h(x)-h(y))$. So the same mapping applies to
  the weak limit, which leaves it invariant so the uniform limits for the triples are equivalent. We have successfully
  changed the last variable in the triple. What remained to change are the two $\tilde\nu$ distributed variables to their
  $\bar\nu$ distributed versions. The $\sigma$-algebra $\bar{\cal S}$, generated by factorised functions, is the
  multiplication of the $\sigma$-algebra generated by the rectangles in $\Delta$ in the stable direction, and the
  Borel-algebra in the unstable direction$\pmod 0$. The
  forthcoming limit theorem for $\bar \Upsilon_n$ proves the same for $F\bar {\cal S}$, because the application of $F$
  means the application of $(x,y,\xi)\mapsto(Fx,Fy,\xi-\bar f(x)+\bar f(y))$, and the limit is invariant under this
  action. Since $\bigvee _{n>0}F^n\bar{\cal S}={\cal S} \pmod 0$ it is enough to prove the limit theorem for $\bar
  \Upsilon_n$. 
  
  For to do this we are going to integrate test functions: $w(\bar x,\bar y,\xi)$. We will restrict ourselves to functions
  which are in $\mathcal{L}$ as functions of $x$ and $y$ and are integrable (with respect to the prospective limit) as
  functions of $\xi$, moreover their Fourier transform is compactly supported. By Breiman \cite{Bre} checking convergence
  for these functions proves weak convergence of measures. For simplicity we are going to use the inverse transform:
  $w(\bar x,\bar y,\xi)=\int\hat w(\bar x,\bar y,t) e^{it\xi}dt$.
  \begin{multline*}
    n^\frac d2 \int_{\bar \Delta\times\bar\Delta\times M(f)} w d\bar\Upsilon_n = n^\frac d2 \int w(\bar x,\bar
    F^n\bar x,S_n(\bar x)-k_n) d\bar\nu \\ 
    = n^\frac d2 \int \int_{\widehat{M(f)}}\hat{w}(\bar x,\bar F^n\bar x,t) e^{it(S_n(\bar x)-k_n)}dt\enskip d\bar\nu \\
    = n^\frac d2 \int \rho^{-1}(\bar x) P^n \rho(\bar x)\left( \int_{\widehat{M(f)}} \hat{w}(\bar x,\bar F^n\bar x,t)
      e^{it(S_n(\bar x)-k_n)}dt \right) d\bar\nu \\
    = n^\frac d2 \int \int_{\text{supp}\hat{w}} \rho^{-1}(\bar x) e^{-itk_n} P^n_t \left( \rho(\bar x) \hat
      w(\bar x, \bar F^n\bar x,t)\right) dt \enskip d\bar\nu
  \end{multline*}
  Using lemma \ref{lem:cdf} and theorem \ref{thm:nag} we can substitute $P^n_t \rho \hat w$ by $ \lambda_t^n \rho_t
  \int_{\bar\Delta} \rho \hat w d\bar m$ in the domain $|t|<\delta$ and we get an error term $O(n^\frac d2 \theta^n)$
  inside the integration wrt $\bar\nu$. This involves the error terms of lemma \ref{lem:cdf} and theorem \ref{thm:nag}.
  Since $\int \hat w d\bar \nu$ depends only on $t$ we will use the shorter $\hat w (t)$ form.
  \begin{multline*}
  n^\frac d2 \int_{\bar \Delta\times\bar\Delta\times M(f)} w d\bar\Upsilon_n = \int \rho^{-1}(\bar x)\int_{\left| t
  \right|<\delta\sqrt{n}} \hat{w}\left(\frac t{\sqrt{n}}\right) e^{-it\frac{k_n}{\sqrt{n}}} \lambda_{\frac t{\sqrt{n}}}^n
  \rho_{\frac t{\sqrt{n}}}(\bar x) dt + o(1) d\bar \nu\\
  \rightarrow \int_{\mathbb{R}^d} \int\hat{w}(\bar x,\bar y,0)d\bar\nu e^{-itk} e^{\frac{-\sigma^2 t^2}2} dt \\
  = \frac 1{(2\pi)^d} \int_{M(f)} w(\bar x,\bar y,\xi) d\bar\nu^2\times dl \quad \frac 1 {\det\sigma} \sqrt{2\pi}^d
  e^{-\frac {k^2}{2\sigma^2}}
  \end{multline*}
  In the above limit the order of the error term is meant in ${\cal L}$-norm (cf.\ lemma \ref{lem:cdf} and theorem
  \ref{thm:nag}), this implies that limiting makes the error term vanish (cf.\ definition of ${\cal L}$-norm). The same
  applies for the $\bar x$ dependence of $\rho_{\frac t{\sqrt n}}$. The convergence in $t$ is dominated, since $\exists C
  \quad \forall |t|\leq \delta\sqrt{n} \quad \left|\lambda_{\frac t{\sqrt{n}}}^n\right| \leq e^{-C|t| ^2}$.
\ep \medbreak

\proc{Remark}
  The case of nonminimal functions is obvious from the first argument of the proof. If $f-g=h-h\circ T$ then the limit
  measure for $f$ differs from the limit measure for $g$ by convolving the distribution of $h$ and of $-h$.
\medbreak

\begin{theo}\label{thm:lecs}
  Let $k_n\in M(f)$ be such that $\frac{k_n - na} {\sqrt{n}}\rightarrow k\in\mathbb{R}^d$, and $\kappa_n\in M(f)$ be such
  that $\frac{\kappa_n - na} {\sqrt{n}}\rightarrow \kappa\in\mathbb{R}^d$. Denote the joint distribution of
  $S_n-k_n,S_m-\kappa_m$ by $\upsilon_{n,m}$! If $f$ is minimal and nondegenerate, then \[ \lim_{n,m,n-m\rightarrow\infty}
  n^{\frac d2} m^{\frac d2}\upsilon_{n,m} \rightarrow \frac {e^{\frac {-k^2}{2\sigma^2}}e^{\frac
      {-\kappa^2}{2\sigma^2}}}{\det^2\sigma(2\pi)^d} l\times l. \]
\end{theo}

\proc{Proof} Again as in the previous proof if we consider the joint distribution $\Upsilon_{n,m}$ of the 5-tuple
$(x,T^nx,T^mx,S_n(x)-k_n,S_m(x)-\kappa_m)$, then it is enough to prove, that \[ \lim_{n,m,n-m\rightarrow\infty} n^{\frac
  d2} m^{\frac d2}\bar\Upsilon_{n,m} \rightarrow \frac {e^{\frac {-k^2}{2\sigma^2}}e^{\frac
    {-\kappa^2}{2\sigma^2}}}{\det^2\sigma(2\pi)^d} \bar\nu^3\times l^2. \] To prove convergence we are going to integrate
test functions: $w(\bar x,\bar y,\bar z,\xi,\zeta)$. Again as in the previous proof we restrict ourselves to the same
class of functions. We are going to use the inverse transform: $w(\bar x,\bar y,\bar z,\xi,\zeta)=\int\hat w(\bar x,\bar
y,\bar z,t,u) e^{i(t\xi+u\zeta)}dt\enskip du$.
  \begin{gather*}
    n^{\frac d2}m^{\frac d2}\int\limits_{\bar \Delta^3\times M(f)^2} w d\bar\Upsilon_{n,m} = n^\frac d2 m^\frac d2 \int
    \int\limits_{\widehat{M(f)}^2} \rho^{-1} e^{-i(tk_n+u\kappa_n)} P^n_t \left( \rho e^{iuS_m}\hat w\right) dt \enskip du
    \enskip d\bar\nu \\ =n^\frac d2 m^\frac d2 \int\int\limits_{\left| t \right|<\delta} \rho^{-1}(\bar x) e^{-itk_n}
    \lambda_t^n \rho_t\int e^{-iu\kappa_n}\int\limits_{\bar\Delta} \rho e^{iuS_m} \hat w d\bar m \enskip du \enskip dt +
    O(n^\frac d2 \theta^n) d\bar\nu 
  \end{gather*}
  Again the inner integration is invariant under $P$, so 
  \begin{align*}
    \int_{\bar\Delta} \rho e^{iuS_m} \hat w d\bar m &= \int_{\bar\Delta} P^m\rho e^{iuS_m} \hat w d\bar m \\ &=
    \int_{\bar\Delta} P_u^m \rho \hat w d\bar m \\ &= \int_{\bar\Delta} \lambda_u^m \rho_u \int_{\bar\Delta}\rho \hat
    wd\bar m +O(\theta^m) d\bar m
  \end{align*}

From this point the variables can be handled separately and the argument of the previous proof should be repeated twice to
get the statement of this theorem.
\ep \medbreak

\section{Recurrence of planar Lorentz-process}
\label{sec:rec}

\subsection{Semi-dispersing billiards}
In this subsection we summarize some basic properties of semi-dispersing billiards. Our aim is to introduce the most
important concepts and fix the notation. For a more detailed description see the literature, especially \cite{KSSz}.

A billiard is a dynamical system describing the motion of a point particle in a connected, compact domain $Q \subset
\mathbb{T}^d$. The boundary of the domain in assumed to be piecewise $C^3$-smooth. Inside $Q$ the motion is uniform while
the reflection at the boundary $\partial Q$ is elastic. As the absolute value of the velocity is a first integral of
motion, the phase space of the billiard flow is fixed as $M=Q\times S^{d-1}$ -- in other words, every phase point $x$ is
of the form $x=(q,v)$ with $q\in Q$ and $v\in \mathbb{R}^d,\ |v|=1$.  The Liouville probability measure $\mu$ on $M$ is
essentially the product of the Lebesgue measures, i.\ e.\ $d\mu= {\rm const.}\, dq dv$. The resulting dynamical system
$(M, S^{\mathbb{R}} , \mu)$ is the (toric) \emph{billiard flow}.

Let $n(q)$ denote the unit normal vector of a smooth component of the boundary $\partial Q$ at the point $q$, directed
inwards $Q$.  Throughout the paper we restrict our attention to \emph{semi-dispersing billiards}: we require for every
$q\in \partial Q$ the second fundamental form $K(q)$ of the boundary component to be non-negative.

The boundary $\partial Q$ defines a natural cross-section for the billiard flow. Namely consider \[ \partial M = \{ (q,v)
\mid q\in \partial Q,\left< v,n(q)\right> \ge 0 \}.\] This set actually has a natural bundle structure (cf.\
\cite{4geom}). The Poincar\'e section map $T$, also called the \emph{billiard map} is defined as the first return map
on $\partial M$. The invariant measure for the map is denoted by $\mu_1$, and we have $d\mu_1= {\rm const.} \left|\left<
    v,n(q)\right>\right|dqdv$. Throughout the paper we work with this discrete time dynamical system $(\partial
M,T,\mu_1)$. Recall the usual notation: for $(q, v) \in M$ one denotes
$\pi(q, v) = q$ the natural projection.

The \emph{Lorentz process} is the natural $\mathbb{Z}^d$ cover of a toric billiard. More precisely: consider
$\Pi:\mathbb{R}^d \to \mathbb{T}^d$ the factorisation by $\mathbb{Z}^d$. Its fundamental domain $D$ is a $d$-dimensional
cube (semi-open, semi-closed) in $\mathbb{R}^d$, so $\mathbb{R}^d = \cup_{z \in \mathbb{Z}^d} (D+z)$, where $D+z$ is the
translated fundamental domain.

By denoting $\tilde Q = \Pi^{-1} Q$, $\tilde M = \tilde Q \times S^{d-1}$, etc., the Lorentz dynamics is $(\tilde M,
\{\tilde {S}^t\mid t \in \mathbb{R}\}, \tilde \mu)$ and its Poincar\'e section map is ($\partial \tilde M, \tilde T,
\tilde \mu_1)$.  The \emph{free flight function} $\tilde\psi: \partial \tilde M \to \mathbb{R}^d$ is defined as follows:
$\tilde \psi(\tilde x)=\tilde q(T\tilde x)-\tilde q (\tilde x)$. The \emph{discrete free flight function} $\tilde \kappa :
\partial \tilde M \to \mathbb{Z}^d$ is defined as follows: $ \tilde \kappa (\tilde x) = \iota (\tilde T \tilde x) - \iota
(\tilde x)$, where $ \iota (\tilde x) = z $ if $ \tilde x \in Dz$. Observe finally, that $\tilde\psi$ and $\tilde \kappa$
are invariant under the $\mathbb{Z}^d$ action, so there are $\psi$ and $\kappa$ functions defined on $\partial M$, such
that $\tilde\psi= \Pi^*\psi$ and $\tilde \kappa = \Pi^* \kappa$. Actually for our purposes it will be more convenient to
choose the fundamental domain in such a way that $\partial \tilde Q \cap \partial D=\emptyset$. In this way $\kappa$ will
be continuous.

\subsection{Minimality of the free flight function}

Start with a simple observation

\begin{lemma} 
\[\kappa  \sim \psi\]
\end{lemma}
\proc{Proof}
  Fix an arbitrary point $w \in D$. For $x=(q, v) \in \partial M$ define $h(x)=w-q$. if $h(Tx) \in D+z$ for some $z \in
  \mathbb Z^d$, then $\kappa(x)=z$, and, of course, \[ \psi(x) = \kappa(x) +h(x) - h(Tx)\]
\ep 
\begin{theo}
  $\kappa$ is minimal in the class of $\psi$.
\end{theo}
\proc{Proof} Suppose the contrary and denote the minimal function by $\kappa'$! Apply the factorisation by the minimal
lattice: $\kappa_f :\partial M\to \mathbb{Z}^d/M(\kappa)$! Then $\kappa_f\sim \kappa'_f$, and $\kappa'_f$ is the constant function. Denote
by $n$ the cardinality of this abelian group $\mathbb{Z}^d/M(\kappa)$! (We can suppose $n<\infty$.) In this case $\forall
x$ periodic, such that $n|\mathrm{per}(x)=p$ the Birkhoff sum $S_p(\kappa_f)x=0$. The proof of the theorem is based on our forthcoming lemma \ref{lem:5.2}.\@ It is a variant of a statement which was originally applied in 
\cite{BSCH91} to establish the non-singularity of the limiting covariance in the CLT.\@ To contradict the non-minimality we are going to find a periodic point for each sublattice of finite index,
not satisfying the above equation.
  \begin{lemma}\label{lem:5.2}
    For any finite index sublattice $Z\subset\mathbb{Z}^d$ there exists a periodic point $x$ such that the period $p$ is a
    multiple of $\left|\mathbb{Z}^d:Z\right|$ and $\sum_{i=0}^{p-1} \kappa(T^ix) \not\equiv 0 \pmod Z$
  \end{lemma}
  \proc{Proof of lemma} The idea is a suitably adapted, simplified and generalized version of an argument of
  \cite{BSCH91}. The original idea is well explained in \cite{BinSz}. Fix the lattice $Z$, denote the index by $i$, and fix $\Lambda \subset \mathbb{T}^d_0$, the basic product
  set of the Young system of our billiard ($\mu_1(\Lambda) >0$).  Take a billiard in the elongated torus $\mathbb{T}(Z) =
  \mathbb{R}^d /Z$, which is an appropriate projection of our Lorentz process. Consider the images of $\Lambda$ on the
  elongated torus. Take two of them $\Lambda_0$ and $\Lambda_1$. By using the ergodicity of powers of the billiard in
  $\mathbb{T}(Z)$ we see that there exists an $n \in \mathbb{Z}_+$ such that $\Lambda_0 \cap T(Z)^{-ni} \Lambda_1$
  contains a Markov intersection $\Lambda^*$ of positive measure where $T(Z)$ denotes the Poincar\'e section map of the
  billiard on $\mathbb{T}(Z)$. The fact that $\Lambda_0 \cup T(Z)^{-ni} \Lambda_1$ contains a Markov intersection
  $\Lambda^*$ of positive measure requires a proof. This is the only part in our paper where we have to go beyond
  properties (P1-8) of Young systems formulated in subsection 2.1 and to use some more detailed arguments from her
  construction. To make the reading of the main body of this paper easier we will postpone until the Appendix the proof of
  the sublemma formulating this particular statement .
    
    \proc{Sublemma}  For the billiard on $\mathbb {T}(Z)$ there exists an $n \in \mathbb{Z}_+$ such that $\Lambda_0 \cap
    T(Z)^{-ni} \Lambda_1$ contains a Markov intersection $\Lambda^*$ of positive measure.\medbreak
    
    By identifying $\Lambda$ with $\Lambda_0$, $\cap_{l=- \infty}^{\infty} T^{lni} \Lambda^*$ consists of exactly one
    point $x^*$. Clearly $T^ni x^* = x^*$ and, moreover, the claim of the lemma is also evident.  \ep \medbreak To conclude the proof it is sufficient to observe that the relation $\kappa\sim \kappa'$ and the periodicity of $x$ also imply that $\sum_{i=0}^{p-1} \kappa(T^ix) \not\equiv 0 \pmod Z$. Hence the
    theorem.  \ep \medbreak

\subsection{Proof of recurrence}
\label{sec:ut}

In this subsection we want to apply the local limit theorem in order to get the recurrence for the planar Lorentz-process, a result
already proved in \cite{Sch} and in \cite{Conze}.\@ Let the system be a billiard on the 2-dimensional torus, with strictly
convex scatterers, and finite horizon. Such a system is always a Young system. This was proved in \cite{Young}. For the
role of $f$ in the main theorem let we choose $\kappa: X\rightarrow \mathbb{R}^2$ the discrete free flight function. Time
reversion symmetry ensures zero average. We have just proved that $\kappa$ is minimal. Its boundedness is equivalent with the
finite horizon assumption, and the other conditions are trivial. Then theorem \ref{thm:main} ensures that $\nu(S_n\in D) >
\frac Cn$ for some $C>0$. It immediately extends to any fixed domain.
\begin{theo}
  The  planar Lorentz process with finite horizon is almost surely recurrent.
\end{theo}
\proc{Proof} The proof follows the ideas used in \cite{KSz1}. The sequence of events 
\begin{gather*}
  A_n = \left\{ S_n \in D \right\} \intertext{fulfills the condition of Lamperti's Borel-Cantelli \cite{S}:}
  \sum_{k=1}^\infty \nu\{A_k\}=\infty \intertext{is clear by the main theorem} \liminf_{n\rightarrow\infty} \frac
  {\sum\limits _{j,k=1}^n \nu(A_jA_k)} {\left( \sum\limits_{k=1}^n \nu(A_k) \right)^2} <c \intertext{the denominator is
    of order $\log^2 n$, the numerator will be decomposed as follows:} \sum_{j,k=1}^n \nu(A_jA_k) \leq
  \sum_{\min(j,k)<\log n} \nu(A_jA_k) + \sum_{|j-k|<\log n} \nu(A_jA_k) + \sum_{j,k,|j-k|\geq\log n} \nu(A_jA_k).
\end{gather*}
The first sum can be estimated by $2\log n \sum_{k=1}^n m(A_k)$ which is of order $\log^2 n$. The same is true for the
second term as well.  Concerning the third one, by theorem \ref{thm:lecs} we know that the asymptotics of the summand is
proportional to $\frac1{jk}$, so the sum is of order $\log^2 n$. Consequently, by Lamperti's lemma \[ \nu \{A_k \text{
  i.\ o.} \} > \frac 1c. \] Since this event is invariant under the ergodic dynamics, it happens almost surely.  \ep

Finally it is interesting to note that, as observed by Sim\'anyi \cite{Sim 89} the recurrence of the planar Lorentz process is equivalent to saying that the corresponding billiard in the whole plane (with an infinite invariant measure) is ergodic (see also \cite{Pene}).
\section*{Appendix: Proof of sublemma}

The only aim of this appendix to provide the proof of Sublemma.

\proc{Sublemma} For the billiard on $\mathbb {T}(Z)$ there exists an $n \in \mathbb{Z}_+$ such that $\Lambda_0 \cap
T(Z)^{-ni} \Lambda_1$ contains a Markov intersection $\Lambda^*$ of positive measure.\medbreak

\proc{Proof} In order that our ideas be clear with a minimal knowledge of sections 7 and 8 of \cite{Young} we summarize
some facts from this reference.  First, let us note that often it is convenient to use the semi-metric $p$ determined by
the density $\cos \phi dr$. We will write $p(.)$ for the $p$-length of a curve, while $l(.)$ denotes its Euclidean
length. Finally, as before, $d(.,.)$ denotes Euclidean distance. In particular, $\gamma^u_\delta(x)$ will denote that
piece of a $\gamma^u_{loc}$-curve whose endpoints have $p$-distance $\delta$ from its `center' $x$.

Facts:
\begin{enumerate}
\item [(i)] $\delta_1 >0$ is a suitably small number, $\delta =\delta_1^4$ and $\alpha_1 = \alpha^{\frac{1}{4}}$.
\item [(ii)] The product set $\Lambda$ has a sort of center $x_0 \in A_{\delta_0}= \{x \in M |\ \gamma^u_{3\delta_0}(x)
  exists\} \neq \emptyset$. Denote $\Omega = \gamma^u_{3\delta_0}(x_0)$. Moreover, let us fix a small, rectangular shaped
  neighbourhood $U$ of $x_0$ such that $\Lambda \cap U$ itself is a product set with $\mu_1(\Lambda \cap U) > 0.$
\item [(iii)] For the product set $\Lambda$ one has a simply connected, rectangular-shaped region $Q(x_0)$ such that
  $\partial Q(x_0)$ is made up of two $u$-curves and two $s$-curves. The two $u$-curves are roughly $2\delta_0$ in length
  and they are either from $\Gamma^u(x_0)$ or do not meet any element of $\Gamma^u(x_0)$. The two $s$-curves are
  approximately $2\delta$ long and have the same properties wrt $\Gamma^s(x_0)$. $\hat Q(x_0)$ is a proper
  $u$-subrectangle of $Q(x_0)$, i.\ e.\ it shares the $s$-boundaries of $Q(x_0)$ and its $u$-boundaries, which must have
  the same properties as those of $Q(x_0)$, are strictly inside $Q(x_0)$.
\item [(iv)] Denote $\Omega_\infty = \{y \in \Omega |\ {\rm for}\ \ \forall \ n \ge 0 \ \ d(T^ny, S) > \delta_1 \alpha^{n}
  \}$.  There are unions of a finite number of closed connected curves $\omega$ such that $\Omega_n \supset \Omega_{n+1}$
  and $\Omega = \cap_n \Omega_n$. In addition, if $\omega$ is a component of $ \Omega_n$, then $T^n\omega$ is a connected
  smooth curve with $d(T^n \omega, S) \ge \frac{1}{2} \delta_1 \alpha^{n}$, and, in particular, $T^{n+1} \omega $ is also
  a connected smooth curve.
\item [(v)] If for a point $x$ one has $R(x) = n$, then $x$ belongs to an $s$-subrectangle $Q_\omega$ of $Q(x_0)$ (where
  $\omega $ is some component figuring in (iv) ) such that $T^j Q_\omega \cap S = \emptyset $ for every $0 \le j \le n$.
  Also, $10 \delta_0 \le p(T^n \omega) \le 20 \delta_0$ and $T^n \omega$ $u$-crosses $\hat Q(x_0)$ with segments $2
  \delta_0$ in length sticking out on both sides.
\item [(vi)] Finally, for some $R_1 \ge R_0$ large enough it is true that if, for some $n \ge R_1$, a component $\omega$
  of $\Omega_n$ $u$-crosses the middle half of $Q$ under $T^n$ , then the entire $s$-subrectangle of Q associated with
  $\omega$ $u$-crosses $Q$ under $T^n$.
\end{enumerate}

When now turning to the billiard on $\mathbb T(Z)$ we will extend our previous usage of notations: for instance,
$x_0^{(0)}, \dots, x_0^{(Z)}$ will denote the different copies of $x_0$, and similarly $U^{(0)}, \dots, U^{(Z)}$ the
different copies of $U$. $\mu_1(Z, .)$ will denote the invariant probability measure for our `elongated' billiard system.
We note that Young's construction uses powers of $T$ which are multiple of some given natural number. Here, for
simplicity, we take this number to be equal to one and use the ergodicity of $T$.  However, for our billiard it is known
that any power of $T$ is also ergodic so our simplification is by no means a restriction.

In fact, claim (vi) is the main fact necessary for our purposes. Introduce the function 
\begin{gather*}
  w(x) = \chi_{\{p(\gamma^u(x))\ge 10 \delta_0\}}(x) \chi_{\{x \in \Lambda^{(Z)} \cap U^{(Z)}\}}(x). \intertext{By
    ergodicity,} \frac{1}{n} \sum_{k=0}^{n-1} \int \chi_{\{x \in \Lambda^{(0)} \cap U^{(0)} \}}(x) w(T^kx) d\mu_1(Z, x)
  \to \mu(\Lambda^{(0)} \cap U^{(0)}) \bar w
\end{gather*}
where $\bar w = \int w(x) d\mu_1(Z, x) >0.$ Therefore, for some $x \in \Lambda^{(0)} \cap U^{(0)}$ there exist arbitrarily
large indices $k$ such that $T^k x \in \Lambda^{(Z)} \cap U^{(Z)}$ and $p(\gamma^u(T^kx)) \ge 10 \delta_0$. Since $x \in
\Omega_\infty^{(0)} \subset \Omega^{(0)}_k$, by property (vi) we are done.

\end{document}